\documentclass[12pt]{article}
\usepackage{amssymb}
\usepackage{amsfonts}
\usepackage{amsmath}
\usepackage{graphicx}
\usepackage[english]{babel}

\newtheorem{theorem}{Theorem}
\newtheorem{lemma}{Lemma}

\newtheorem{corollary}{Corollary}
\newtheorem{definition}{Definition}
\newtheorem{proposition}{Proposition}
\newtheorem{example}{Example}
\begin{document}

\def\LN{\operatorname{LN}}
\def\LW{\operatorname{LW}}
\def\Exp{\operatorname{Exp}}
\def\Weibull{\operatorname{Weibull}}
\def\Pareto{\operatorname{Pareto}}
\def\FbarO{\overline{F_0\mbox{\hskip-1.2pt}}\mbox{\hskip1.2pt}}
\def\Fbarl{\overline{F_1\mbox{\hskip-1.8pt}}\mbox{\hskip1.8pt}}



\title{On discrimination between classes\\ of distribution tails}



\author{I. V. Rodionov\thanks{Faculty of Mechanics and Mathematics, Lomonosov Moscow State University,
Moscow, Russia; Department of Innovation and High Technology, Moscow Institute of
Physics and Technology (State University), Moscow, Russia;
e-mail: {vecsell@gmail.com}}
}
\date{}

\maketitle

\begin{abstract}
We propose the test for distinguishing between two classes of
distribution tails using only the largest order statistics of the sample
and state its consistency. We do not assume belonging the
corresponding distribution functions to any maximum domain of
attraction.
\end{abstract}

\section{Introduction}

Let $X_1,\ldots,X_n$ be independent identically distributed (i.i.d.)
random variables with the continuous distribution function (d.f.)
$F$. We set $x^*_F=\sup\{x:\: F(x)<1\}$. Let $x^*_F=+\infty$. We say
that the tail of the d.f. $H$ is lighter than the tail of the d.f.
$G$ with $x^*_H=x^*_G= +\infty$, if the following condition holds
\[
\frac{1-H(x)}{1-G(x)}\to 0\quad\text{}\ x\to+\infty.
\]

This paper is concerned with the problem of distinguishing between
two arbitrary classes of distribution tails $A_0$ and $A_1$, where
the tails of the distributions lying in one class are lighter than
the tails of the distributions lying in another one. The test of
discrimination between classes of distribution tails proposed in
this paper is asymptotic, we also state its consistency. For the aim
of developing the discrimination test we consider the auxiliary
problem of distinguishing between the simple hypothesis about the
distribution tail and two composite alternatives that include almost
all distributions with the tails lighter or heavier than the
distribution tail of the null hypothesis. We emphasize that unlike
the overwhelming number of works concerned with the problems of
testing hypotheses about the distribution tails we do not assume
that the distribution of the sample should satisfy the conditions of
the Fisher-Tippett-Gnedenko limit theorem, i.e. belong to any
maximum domain of attraction (see the definitions below).

In statistics, one often encounters the problem of discrimination
between close distributions from truncated or censored data -- in
particular, in fields related to insurance, reliability,
telecommunications, computer science and earth sciences. The problem
when only the observations over some threshold are known is well
studied (see the works \cite{Dufour,Guilbaud,Chernobai} and
references therein and the book \cite{Beirlant2}). On the other
hand, according to the statistics of extremes (see the book
\cite{Dehaan}), only higher order statistics can be used for the
problem of discrimination of close tails of distributions, whereas
moderate sample values can be simulated using standard statistical
tools.

Gnedenko's limit theorem (or the extreme value theorem) (see
\cite{Gnedenko}), which is the central result in the stochastic
extreme value theory, states that if there exist sequences of
constants $a_n>0$ and $b_n$ such that the d.f. of normalized maximum
$M_n=\max(X_1,\ldots,X_n)$ tends to some non-degenerate d.f. $G$,
i.e.,
\begin{equation}\label{Gned}
\lim\limits_{n\to\infty} P(M_n\le a_n x+b_n)=G(x),
\end{equation}
then there exist constants $a>0$ and $b$ such that
$G(ax+b)=G_\gamma(x)$, where
\[
G_\gamma(x)=\exp \bigl(-(1+\gamma x)^{-1/\gamma}\bigr),\quad 1+\gamma x>0,\quad
\gamma\in\mathbb{R},
\]
and for $\gamma=0$ the right-hand side should be understood as
$\exp(-e^{-x})$. The parameter $\gamma$ is called the extreme value
index \cite{Dehaan}. The d.f. of a sample $(X_1,\ldots,X_n)$ is said
to belong to the Fr\'echet (Weibull, Gumbel, respectively) maximum
domain of attraction if \eqref{Gned} holds for $\gamma>0$
($\gamma<0$, $\gamma=0$, respectively). The distribution functions
belonging to the Fr\'echet and Gumbel maximum domain of attraction
are called the heavy-tailed and light-tailed distributions
respectively. The distributions with the tails heavier than the
tails of the distributions belonging to the Fr\'echet maximum domain
of attraction are called the distributions with super-heavy tails;
these distributions do not belong to any maximum domain of
attraction (see for details \cite{Bingham}). For investigation the
rates of convergence in the Gnedenko's limit theorem and solving
some other problems of extreme value theory the second extreme value
index $\rho$ is considered (see \cite{HaanResnick}), detailed
investigation of which is beyond the scope of our paper.

The estimators of the extreme value indices $\gamma$ and $\rho$ (see
for details \cite{Dehaan}) can be used in the problem of
discrimination between close distribution tails. In this connection
we refer to \cite{Pickands,Hall,Beirlant,Martins,FHL,Drees,Fraga1},
among many others. Another approach is to estimate the distribution
tails directly using higher order statistics (see
\cite{Smith,HaanSinha,DreesHaan}).

It is clear that the asymptotically normal estimators of the extreme
value index $\gamma$ can be used in the problem of discrimination
between the tails of the distributions belonging to the Fr\'echet
and Weibull domains of attraction respectively. However, the above
approach does not work for a huge class of distribution, for
instance, belonging to the Gumbel maximum domain of attraction
(because in this case $\gamma=0$) or super-heavy-tailed
distributions (because $\gamma$ is not determined). In this
connection we mention the work \cite{Fraga}, the authors of which
propose the test distinguishing between the heavy-tailed and
super-heavy-tailed distributions.

The Weibull and log-Weibull classes of distributions form an
important class of distributions from the Gumbel maximum domain of
attraction. We say that the d.f. $F$ is a Weibull-type d.f., if
there exist $\theta>0$ such that for all $\lambda>0$
\[
\lim\limits_{x\to\infty} \frac{\ln (1-F(\lambda x))}{\ln (1-F(x))}=\lambda^{\theta}.
\]
The parameter $\theta$ is also called the Weibull tail index. The
class of Weibull distributions contains, in particular, the normal,
exponential, gamma distributions and other ones of great value in
statistics. If for some d.f. $F$ the distribution function $F(e^x)$
belongs to the Weibull class with $\theta>1$, then one says that $F$
belongs to the log-Weibull class of distributions. A method capable
of discriminating between close tails of Weibull and log-Weibull
distributions was proposed in \cite{Girard} and based on the well-known
Hill estimator (see \cite{Hill}, and Section~2 of this work), the
estimators of the Weibull tail index (see \cite{Berred,Girard2}) can
also be applied for discriminating between tails of Weibull-type
distributions. The likelihood ratio method applied to higher order
statistics of the sample was used in \cite{Rodionov1,Rodionov2} to
develop the tests discriminating between the tails of the
distributions from the Gumbel maximum domain of attraction.

This paper is actually the first attempt to develop the general test
in the problem of discrimination between the distribution tails that
are not assumed to belong to any maximum domain of attraction. The
problem of optimality of this test is natural and will be considered
in our next works. The problem of the optimal choice of the number
of retained higher order statistics is considered in Section~3.

The paper is organized as follows. In Section~2 we formulate the
problem and propose the test of discrimination. In Section~3 we
illustrate the numerical performances of the test and compare it
with some other tests. The results and the proofs are given in
Sections~4,~5.

\section{The model and the test of discrimination}

Let as above $X_1,\ldots,X_n$ be i.i.d. random variables with a
common continuous d.f. $F$ and $x_F^*=+\infty$. Let
$X_{(1)}\le\ldots\le X_{(n)}$ be the order statistics pertaining to
$X_1,\ldots,X_n$. For the purpose of developing the discrimination
test consider the Hill-type statistic
\[
R_{k,n}=\ln(1-F_0(X_{(n-k)}))-\frac{1}{k}\sum\limits_{i=n-k+1}^n
\ln(1-F_0(X_{(i)})),
\]
where $F_0$ is a continuous d.f. Note that if
$F_0(x)=(1-\kappa^{1/\gamma}/ x^{1/\gamma}) I(x>\kappa)$ is the
Pareto d.f., then
\[
R_{k,n}\overset{d}{=}\gamma_{H}/\gamma
\]
as $X_{(n-k)}>\kappa$, where $\gamma_{H}$ is the Hill's estimator of
the extreme value index $\gamma$ (which in this case is equal to the
parameter of the Pareto distribution, \cite{Dehaan}),
\[
\gamma_H=\frac{1}{k}\sum\limits_{i=n-k+1}^n \ln X_{(i)}-\ln X_{(n-k)},
\]
the estimator is consistent for positive values of $\gamma$.

\begin{definition}\label{def1}
We say, that distribution functions $H$ and $G$ satisfy the
condition $B(H,G)$ (written B-condition), if for some
$\varepsilon\in(0,1)$ and $x_0$
\begin{equation}\label{cond}
\frac{(1-H(x))^{1-\varepsilon}}{1-G(x)}\ \text{ is nonincreasing as
}\ x>x_0.
\end{equation}
\end{definition}

It is easy to see, that under this condition the tail of the
d.f.~$H$ is lighter than the tail of the d.f. $G$, i.e.
\[
\frac{1-H(x)}{1-G(x)}\to 0\quad\text{as}\ x\to\infty.
\]
For example, if $H(x)=1-\exp\{-\lambda x\}$, $x>0$, and
$G(x)=1-\exp\{-\theta x\}$, $x>0$, are two exponential distribution
functions with $\lambda> \theta$, then the concerned condition is
satisfied for all $\varepsilon\le 1- \theta/\lambda$ and $x_0>0$.
The B-condition is also satisfied for the normal distribution
functions with different variances. In addition note that the tail
of the log-Weibull type distributions are lighter then the Pareto
tails but heavier then the tails of the Weibull-type distributions
(see \cite{Girard}), and the B-condition holds for all these classes
of distributions.

Let us formulate the main problem of our paper. Consider two
distribution classes $A_0$ and $A_1$ such that the tails of the
distributions lying in $A_0$ are lighter than the distribution tails
lying in $A_1$. Assume, that there exist the ``separating'' d.f.
$F_0$ such that both of the conditions $B(G, F_0)$ and $B(F_0, H)$
are satisfied with some (may be, different) $\varepsilon$ and $x_0$
for all $G\in A_0$ and $H\in A_1$. Consider in details the examples
of fulfilment of these conditions.

\begin{example}
Let the common distribution of the sample belong to some
one-parametric distribution family, $F\in\{G_\theta,\:
\theta\in\Theta\}$, and for all $\theta_1,\theta_2\in\Theta$,
$\theta_1<\theta_2$, the distribution functions $G_{\theta_1}$ and
$G_{\theta_2}$ satisfy the condition $B(G_{\theta_1},
G_{\theta_2})$. Then, taking $F_0=G_{\theta_0}$, we can distinguish
between the classes $A_0=\{G_\theta,\: \theta< \theta_0\}$ and
$A_1=\{G_\theta,\: \theta> \theta_0\}$.
\end{example}

\begin{example}
We can select $F_0(x)=(1-\exp\{-\exp\{(\ln x)^{1/2}\}\}) I(x>1)$ as
the separating function in the problem of discrimination between
distributions with the Weibull and log-Weibull tails. Indeed, we can
represent the Weibull-type d.f. in the form
\[
F_W(x)=1-\exp\{-\exp\{\theta \ln x (1+o(1))\}\}\quad \text{as}\
x\to\infty,
\]
whereas the log-Weibull-type d.f. can be represented as
\[
F_{\LW}(x)=1-\exp\{-\exp\{\theta \ln\ln x (1+o(1))\}\}\quad \text{as}\ x\to\infty,
\]
hence the fulfilment of the B-conditions is required.
\end{example}

\begin{example}
In the problem of discrimination between heavy-tailed distributions
(i.e. belonging to Fr\'echet maximum domain of attraction) and
super-heavy-tailed distributions (i.e. with the tail heavier than
the tail of any heavy-tailed distribution), the separating function
cannot be correctly selected, but the problem can be solved if we
consider a more particular problem. According to \cite{Dehaan}, the
tail of the arbitrary d.f. belonging to the Fr\'echet maximum domain
of attraction can be represented as $1-F(x)=x^{-1/\gamma} R(x)$,
where $R(x)$ is a slowly varying function, i.e. such that
$\lim\limits_{x\to\infty} R(tx)/R(x)=1$, and the parameter $\gamma$
coincides with the extreme value index of this distribution. So if
one consider the problem of discrimination between
super-heavy-tailed distributions and distributions with heavy tails,
the extreme value index of which is less than some parameter
$\gamma_0$, then the d.f. $F_0(x)=1-x^{-1/2\gamma_0}$, $x>1,$ can be
selected as a separating function.
\end{example}

Thus, suppose that two classes $A_0$ $A_1$ of distributions are
separable by the function $F_0$. Consider the null hypothesis
$\tilde{H}_0\colon F\in A_0$ and the alternative $\tilde{H}_1\colon
F\in A_1$. Propose the following procedure testing the hypothesis
$\tilde{H}_0$:
\begin{equation}\label{test}
\text{}\ R_{k,n}>1+ \frac{u_{1-\alpha}}{\sqrt{k}},\: \text{}\ \tilde{H}_0\ \text{}.
\end{equation}
According to the corollary of Theorem \ref{Th2} (see below), the
proposed test is consistent against the alternative $\tilde{H}_1$
and has the asymptotic significance level $\alpha$, here
$u_{1-\alpha}$~is the quantile of the corresponding level of a
standard normal distribution.

If the tails of distributions lying in $A_0$ are heavier than the
tails of distributions from $A_1$, and there exist the d.f. $F_0$
such that both the conditions $B(F_0, G)$ and $B(H, F_0)$ are
satisfied for some (may be different) $\varepsilon$, $x_0$ and
arbitrary distribution functions $G\in A_0$ and $H\in A_1$, then the
test discriminating between the hypothesis $\tilde{H}_0\colon F\in
A_0$ and the alternative $\tilde{H}_1\colon F\in A_1$ with the same
asymptotic properties is the following:
\begin{equation}\label{test1}
\text{}\ R_{k,n}<1+ \frac{u_{\alpha}}{\sqrt{k}},\: \text{}\ \tilde{H}_0\ \text{}.
\end{equation}

\section{Simulation study}

The aim of the present simulation study is to illustrate the use of
the proposed test and to demonstrate its asymptotic properties.
Firstly consider the problem of discrimination between the Weibull
and log-Weibull classes of distribution tails applying the proposed
test \eqref{test}, and we select the following separating function
$F_0(x)=\bigl(1-\exp(-\exp(\sqrt{\ln x}))\bigr) I(x>1)$. As it
follows from Theorem ~2, the statistics $\sqrt{k}(R_{k,n}-1)$
converges in probability to $-\infty$ on distributions belonging to
the Weibull class and to $+\infty$ on distributions belonging to the
log-Weibull class, that is confirmed by simulations (see Fig.~1).

\begin{figure}[h]
\begin{minipage}[h]{0.49\linewidth}
\center{\includegraphics[width=1\linewidth]{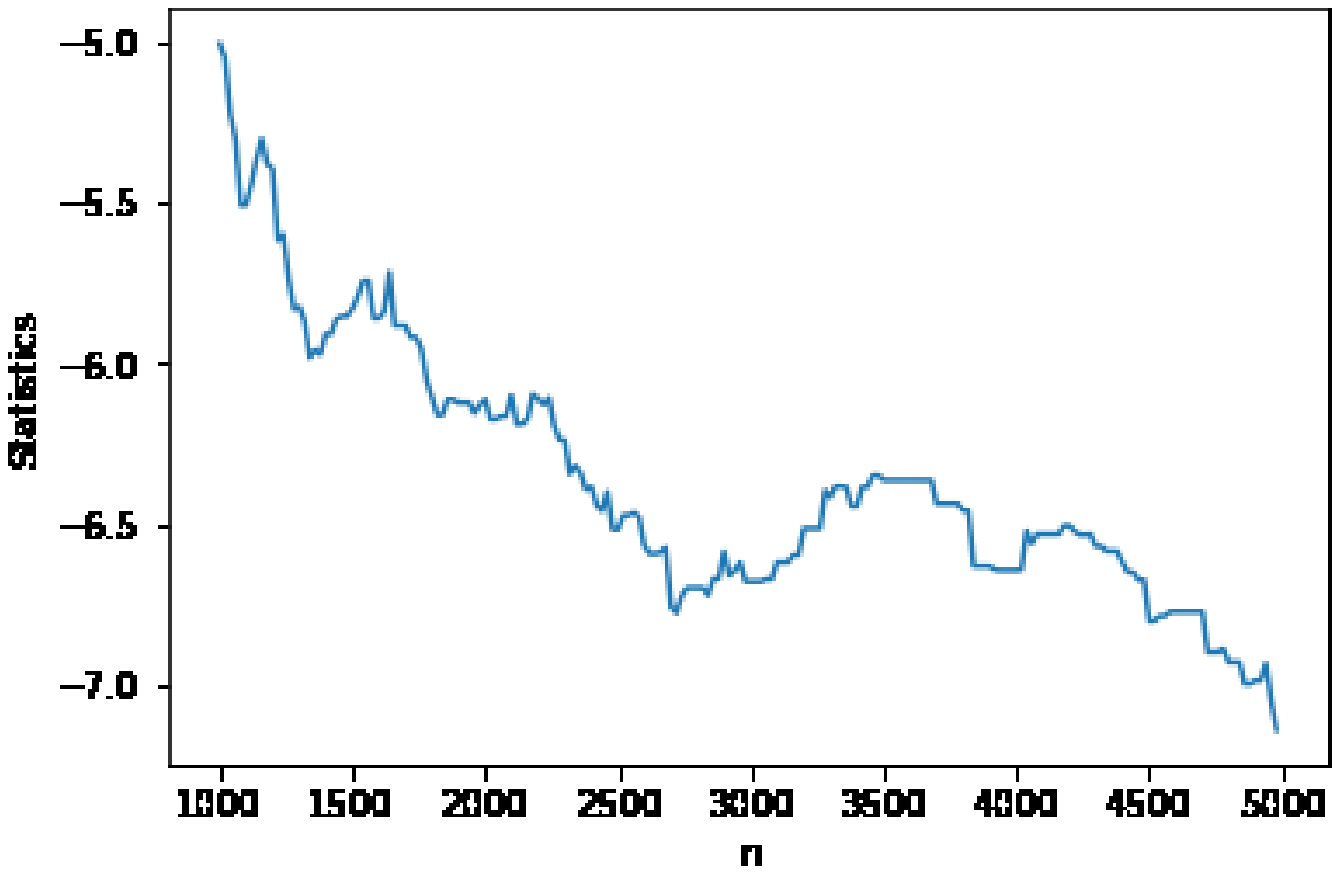} \\ N(0,1)}
\end{minipage}
\hfill
\begin{minipage}[h]{0.49\linewidth}
\center{\includegraphics[width=1\linewidth]{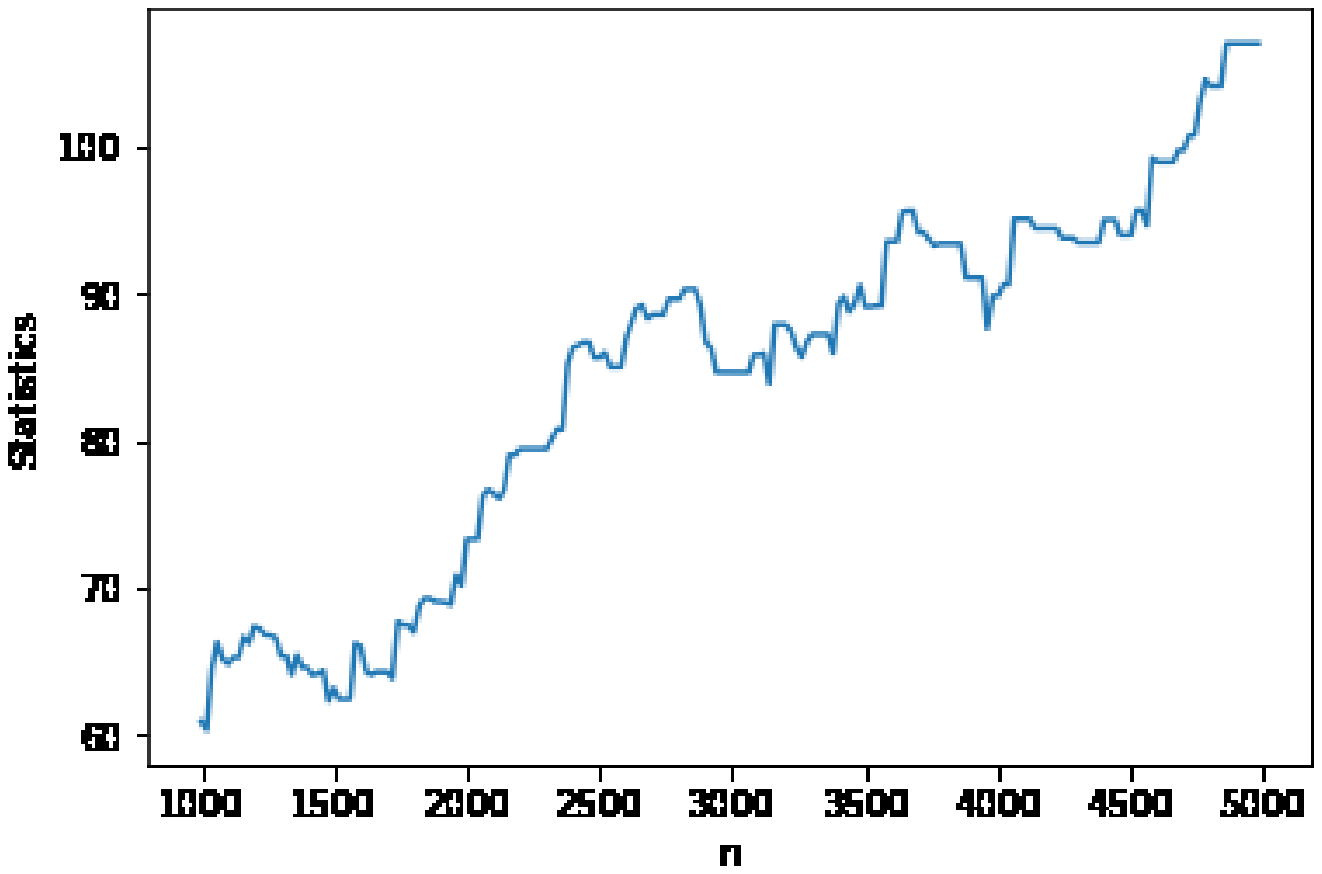} \\
LN(0,1)}
\end{minipage}
\caption{\footnotesize Fig. 1. Values of the statistic $\sqrt{k}(R_{k,n}
- 1)$ for samples obeying $N(0,1)$ distribution (left) and $LN(0,1)$ distribution (right), $n\in [1000; 5000],$ $k=100.$}
\end{figure}

Empirical type I error probabilities and empirical power of the test
\eqref{test} discriminating between the Weibull and log-Weibull
classes for different distributions with the separating function
$F_0(x)=\bigl(1-\exp(-\exp(\kern-1pt\sqrt{\ln x}))\bigr)I(x>1)$ are
given in Table~1 (the nominal level $\alpha=0.05$). We include the
simulation results for the log-Weibull distribution $\LW(3,1)$ with
the d.f. $F_{\LW}(x)=(1-\exp(-(\ln x)^3))I(x>1)$, the log-normal
distribution with parameters $(0,1)$ and the standard exponential
distribution. For the standard normal distribution and the
distributions with lighter tails the empirical type I error
probabilities are significantly less than $0.01$ for all considered
values of $n$ and $k$, so we do not include these probabilities in
the Table~1. In addition, as an example we include the empirical
power of the test for Pareto distribution with the parameter 2 and
standard Cauchy distribution.

\begin{table}[tp]
\begin{tabular}{|l|c|c|c|c|c|c|}
\multicolumn{7}{l}{\parbox{.97\textwidth}{\textbf{Table 1. }
Empirical type I error probabilities an empirical power of the test
\eqref{test} for various distributions, built on
$m=10000$ samples.}}\\[2.5ex] \hline & $n=100$\rule{0pt}{10pt} & $n=100$ & $n=200$ & $n=200$
& $n=500$ & $n=500$\\ & $k=10\phantom{1}$ & $k=20\phantom{1}$ & $k=20\phantom{1}$ &
$k=50\phantom{1}$ & $k=50\phantom{1}$ & $k=100$\\ \hline $\LW(3,1)$\rule{0pt}{10pt} &
$0.29$ & $0.43$ & $0.48$ & $0.75$ & $0.8\phantom{1}$ & $0.95$\\ $\LN(0,1)$ & $0.66$ &
$0.8\phantom{1}$ & $0.89$ & $0.97$ & $0.99$ & $1$\\ $\Exp(1)$ & $0$ & $0.01$ & $0$ &
$0.01$ & $0$ & $0$\\ $\Pareto(2)$ & $0.63$ & $0.61$ & $0.86$ & $0.79$ & $0.99$ &
$0.99$\\ Cauchy & $0.99$ & $0.99$ & $1$ & $1$ & $1$ & $1$\\ \hline
\end{tabular}
\end{table}

Now consider the problem of discrimination between the hypothesis
$H_0\colon \theta<2$ and the alternative $H_1\colon \theta>2$, where
$\theta$ is the Weibull-tail index. Let us compare the test
\eqref{test} (we select as a separating function $F_0(x)$ the
standard normal d.f.), and the test proposed in work \cite{Girard}
(see in addition \cite{GG}). If we select the sample from the
standard normal distribution, then the limit distributions of both
the test statistic \eqref{test} and the test statistic proposed in
\cite{Girard} is standard normal (see the second plot in Fig.~2).

\begin{figure}[h]
\begin{minipage}[h]{0.32\linewidth}
\center{\includegraphics[width=1\linewidth]{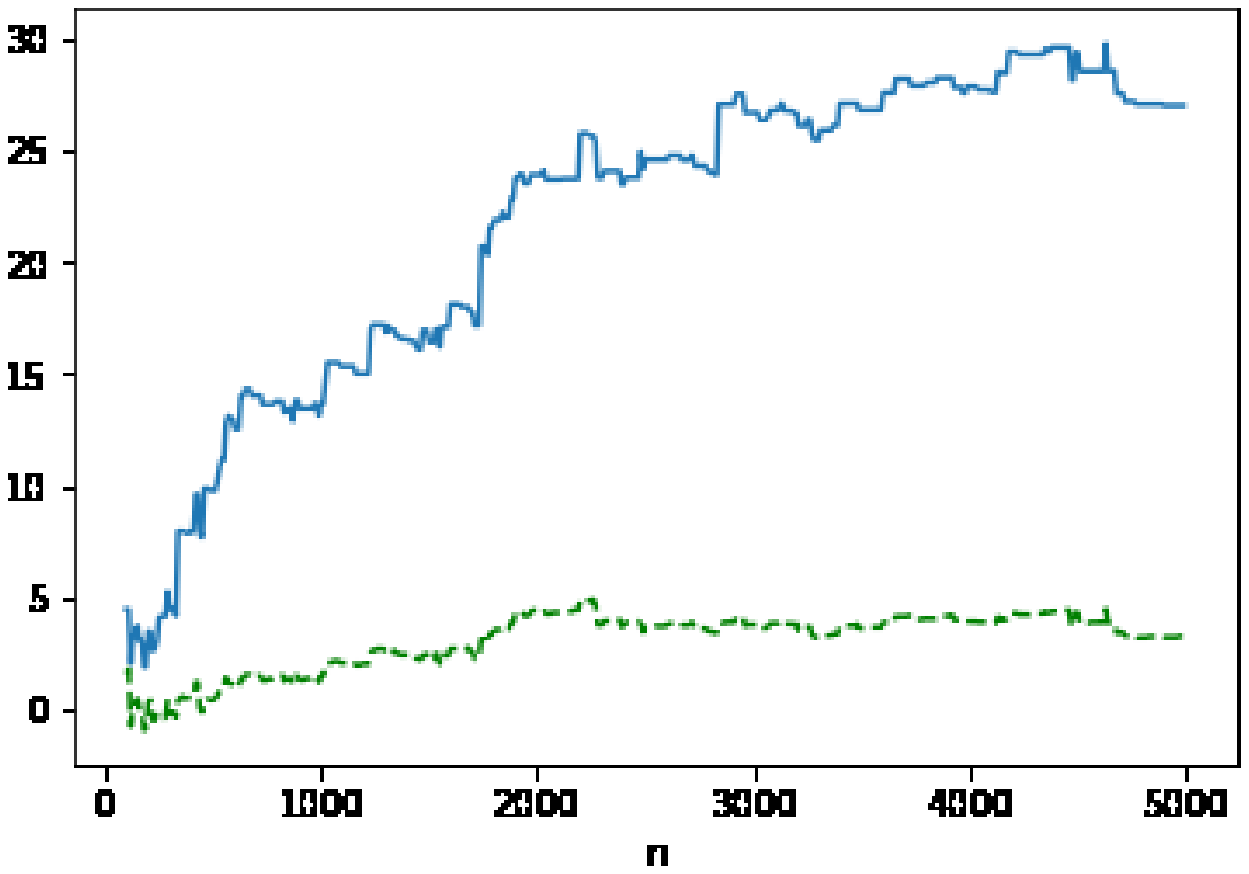} \\
Exp(1)}
\end{minipage}
\hfill
\begin{minipage}[h]{0.32\linewidth}
\center{\includegraphics[width=1\linewidth]{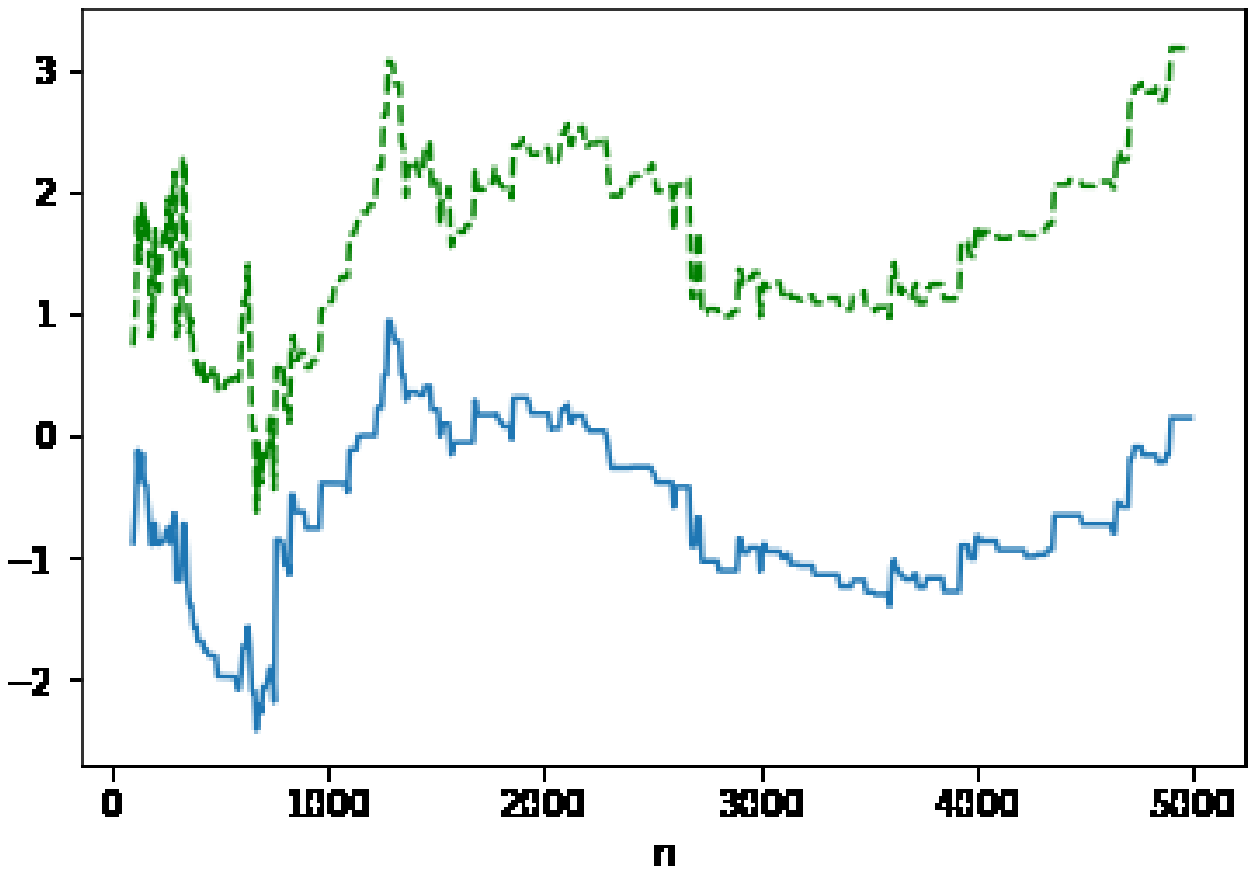} \\
N(0,1)}
\end{minipage}
\hfill
\begin{minipage}[h]{0.32\linewidth}
\center{\includegraphics[width=1\linewidth]{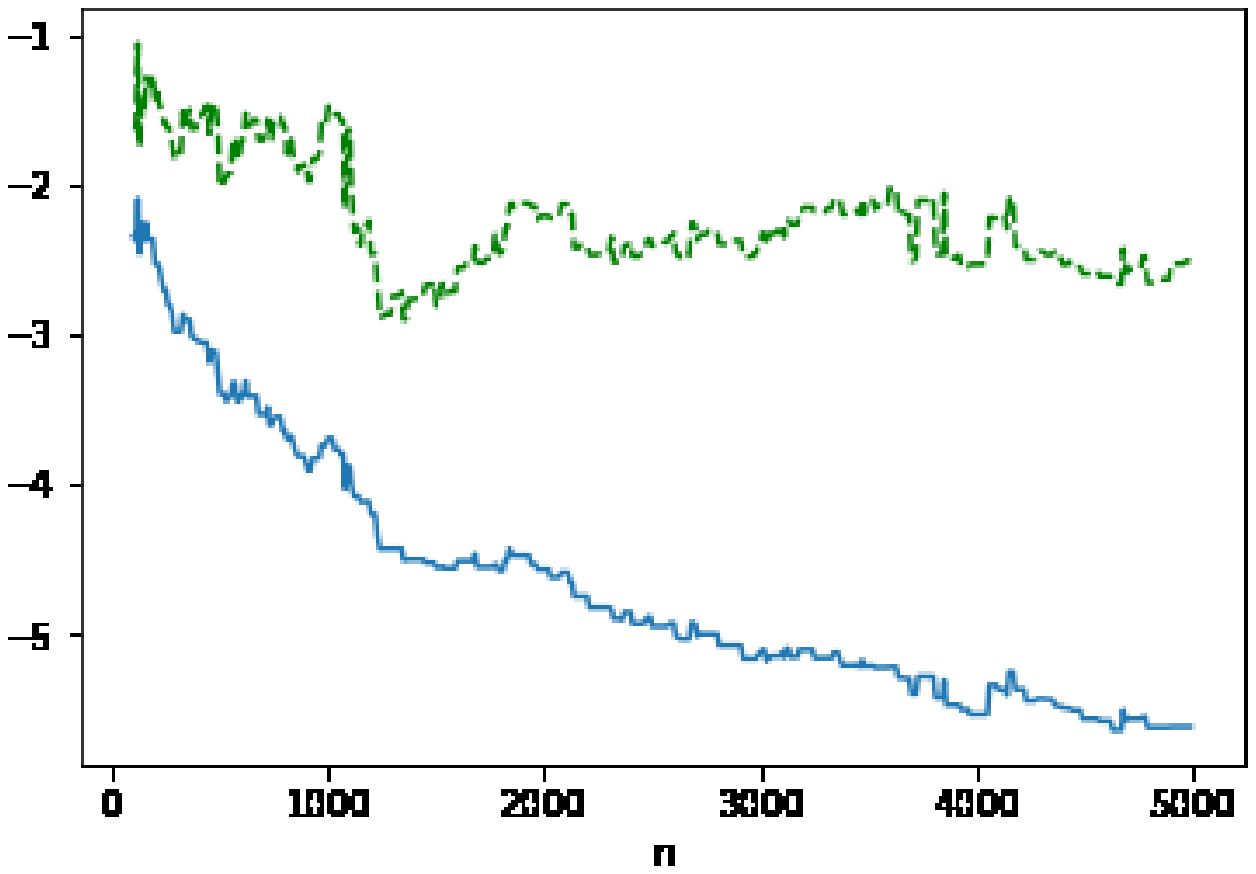} \\
Weibull(3,1)}
\end{minipage}
\caption{\footnotesize Fig. 2. Values of statistics of the test
\eqref{test} (solid line) and test proposed in \cite{GG} (dash-dotted line) on the distributions $Exp(1),$ $N(0,1),$ $Weibull(3,1),$ $n\in[100, 5000],$ $k = [5\ln n].$}
\end{figure}

Consider the problem of distinguishing between distributions with
heavy and super-heavy tails. As it is mentioned before, the proposed
test cannot be applied in this problem because of impossibility of
selecting a separating function $F_0$. But for more particular
problem of discriminating between distributions with heavy tails and
logarithmic distributions with the d.f. $F(x)=1-(\ln
x)^{-\theta}(1+o(1))$, $\theta>0$, as $x\to\infty$, the proposed
test is applicable. Let $A_0$ be the class of logarithmic
distributions mentioned above and $A_1$ be the class of heavy-tailed
distributions. Consider the problem of discrimination between the
hypotheses $H_0\colon F\in A_0$ and $H_1\colon F\in A_1$. We select
as a separating function for the test \eqref{test1} the following
function
\[
F_0(x)=1-\exp(-\sqrt{\ln x}).
\]
In Tables 2 and 3 one can find the empirical type I error
probabilities and empirical power of the test \eqref{test1} and test
proposed in \cite{Fraga} respectively, for the log-Pareto
distribution identified by the d.f. $F_{LP}(x)=(1-(\ln
x)^{-\theta})I(x>e)$ with the parameters $\theta=1$ and $\theta=2$,
the log-$\operatorname{Gamma}(2,1)$ distribution with the density
function $f_{LG}(x)=x^{-2}\ln x I(x>1)$, the standard Cauchy
distribution, the Pareto distribution with the parameter 2 and
standard log-normal distribution. Note the high values of the type I
error probabilities of the test proposed in \cite{Fraga} for the
log-$\Pareto(2)$ distribution.

\begin{table}[tp]
\begin{tabular}{|l|c|c|c|c|c|c|c|}
\multicolumn{8}{l}{\parbox{0.97\textwidth}{\textbf{Table 2. }
Empirical type I error probabilities and empirical power of the test
\eqref{test1} for various distributions, built on
$m=10000$ samples.}}\\[2.5ex] \hline & $n=100$\rule{0pt}{10pt} & $n=100$ & $n=200$ & $n=200$
& $n=500$ & $n=1000$ & $n=5000$\\ & $k=10$ & $k=20$ & $k=20$ & $k=50$ & $k=50$ &
$k=50$ & $k=50$\\ \hline log-$\Pareto(2)$\rule{0pt}{10pt} & $0.07$ & $0.2\phantom{1}$
& $0.06$ & $0.36$ & $0.03$ & $0$ & $0$\\ log-$\Pareto(1)$ & $0$ & $0$ & $0$ & $0$ &
$0$ & $0$ & $0$\\ log-Gamma & $0.56$ & $0.78$ & $0.9\phantom{1}$ & $0.98$ & $1$ & $1$
& $1$\\ Cauchy & $0.62$ & $0.88$ & $0.94$ & $1$ & $1$ & $1$ & $1$\\ $\Pareto(2)$ &
$0.98$ & $1$ & $1$ & $1$ & $1$ & $1$ & $1$\\ $\LN(0,1)$ & $0.99$ & $1$ & $1$ & $1$ &
$1$ & $1$ & $1$\\ \hline
\end{tabular}
\vskip\bigskipamount
\begin{tabular}{|l|c|c|c|c|c|c|c|}
\multicolumn{8}{l}{\parbox{0.97\textwidth}{\textbf{Table 3. }
Empirical type I error probabilities and empirical power of the test
proposed in \cite{Fraga} for various distributions, built on
$m=10000$ samples.}}\\[2.5ex] \hline & $n=100$\rule{0pt}{10pt} & $n=100$ & $n=200$ & $n=200$
& $n=500$ & $n=1000$ & $n=5000$\\ & $k=10$ & $k=20$ & $k=20$ & $k=50$ & $k=50$ &
$k=50$ & $k=50$\\ \hline log-$\Pareto(2)$\rule{0pt}{10pt} & $0.28$ & $0.71$ & $0.49$
& $0.98$ & $0.84$ & $0.6\phantom{1}$ & $0.19$\\ log-$\Pareto(1)$ & $0.06$ & $0.12$ &
$0.07$ & $0.3\phantom{1}$ & $0.09$ & $0.06$ & $0.04$\\ log-Gamma & $0.37$ & $0.58$ &
$0.62$ & $0.89$ & $0.94$ & $0.97$ & $0.99$\\ Cauchy & $0.51$ & $0.8\phantom{1}$ &
$0.8\phantom{1}$ & $0.99$ & $0.99$ & $1$ & $1$\\ $\Pareto(2)$ & $0.92$ & $1$ & $1$ &
$1$ & $1$ & $1$ & $1$\\ $\LN(0,1)$ & $0.91$ & $0.98$ & $1$ & $1$ & $1$ & $1$ & $1$\\
\hline
\end{tabular}
\end{table}

The problem of selecting the number of retained order statistics for
constructing the estimators and tests is significant, but one of the
most complicated in statistics of extremes (see the book
\cite{Fraga2} and references therein). In the context of the problem
of distinguishing between heavy-tailed and super-heavy-tailed
distributions consider the behavior of the empirical type I error
probability and empirical power of the test \eqref{test1} with the
separating function $F_0(x)=1-\exp(-\sqrt{\ln x})$ as $k$ increases
(see Fig.~3). We see, that with the level 0.05 the optimal value of
$k$ lies between 50 and 150.

\begin{figure}[h]
\begin{minipage}[h]{0.49\linewidth}
\center{\includegraphics[width=1\linewidth]{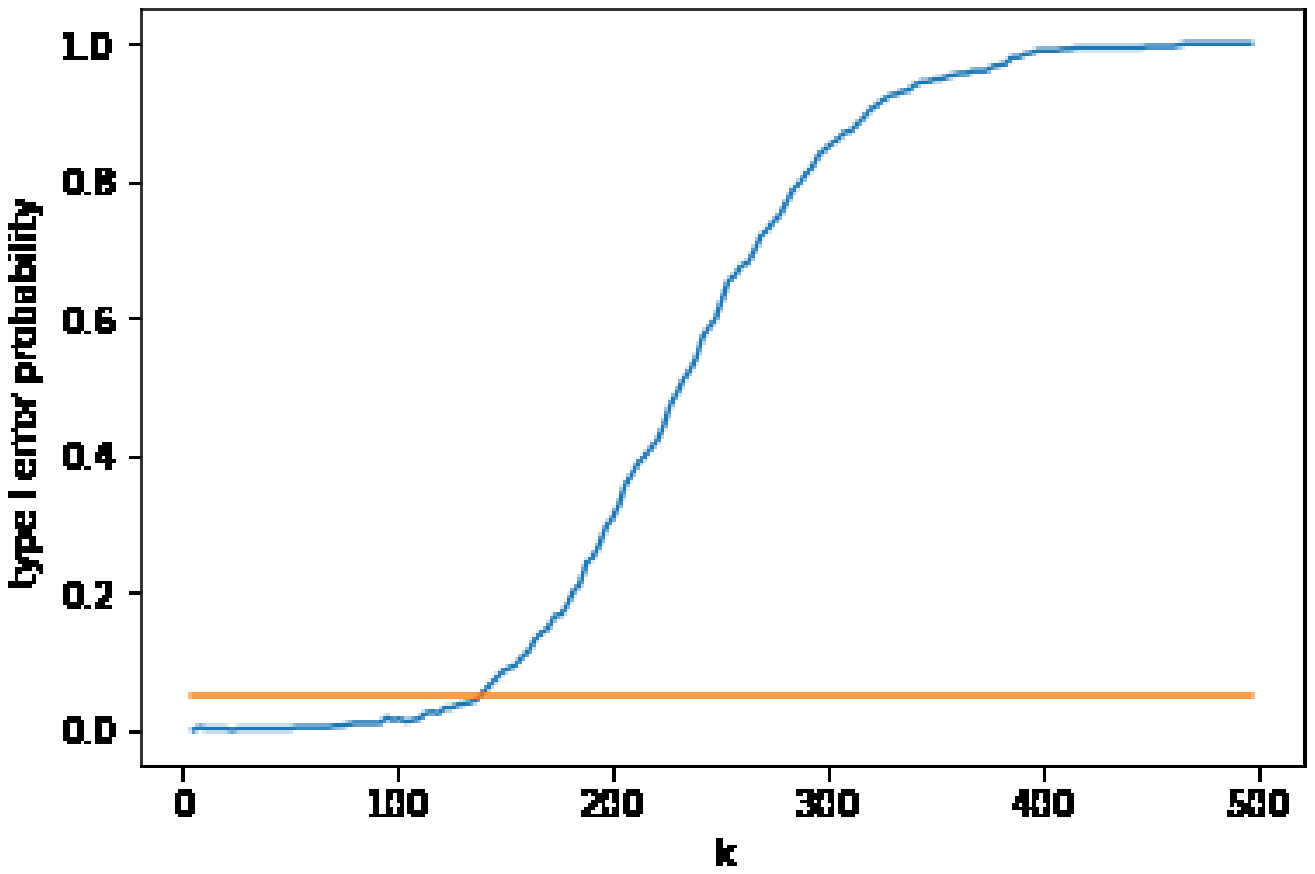} \\
log-Pareto(2)}
\end{minipage}
\hfill
\begin{minipage}[h]{0.49\linewidth}
\center{\includegraphics[width=1\linewidth]{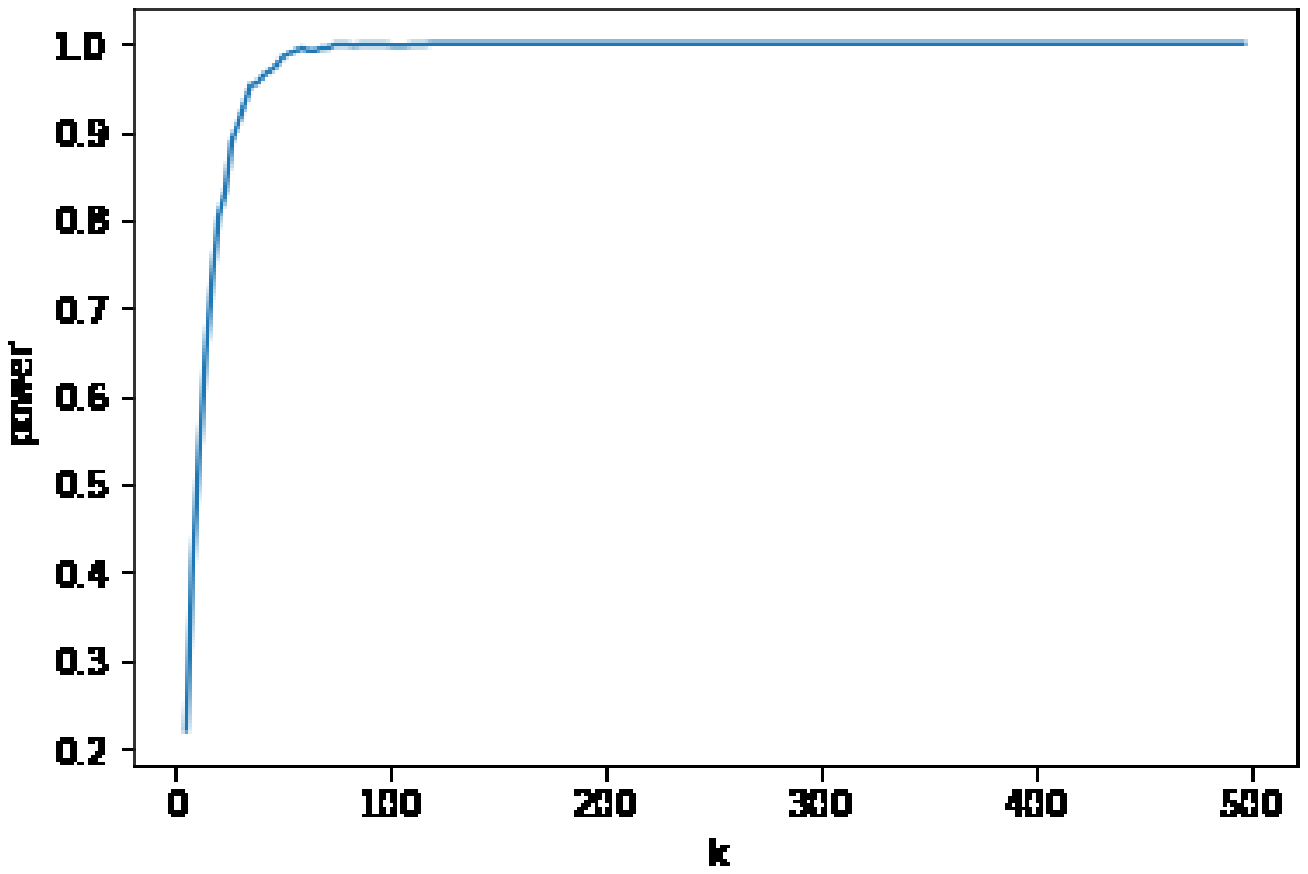} \\
Pareto(0.5)}
\end{minipage}
\caption{\footnotesize Fig. 3. Empirical type 1 error probabilities (left) and empirical power (right) of the test \eqref{test1} on the distributions log-Pareto(2) and Pareto(0.5), respectively, $n=1000.$}
\end{figure}

\section{Main results}

This section discusses theoretical properties of the test proposed
in Section~2. As before, $X_1,\ldots,X_n$ are i.i.d. random
variables with the continuous d.f. $F$. Let $\Theta(F_0)$ be the
class of the continuous distribution functions $F_1$ satisfying
either $B(F_1,F_0)$ or $B(F_0,F_1)$. Consider the simple hypothesis
$H_0\colon F=F_0$ (actually we check the hypothesis that $F$ and
$F_0$ have the same tail) and the alternative $H_1\colon
F\in\Theta(F_0)$. Note that if $F_0$, $F_1$ satisfy either
$B(F_1,F_0)$ or $B(F_0,F_1)$ for some $\varepsilon>0$, then the same
holds for all $\varepsilon_1$, $0<\varepsilon_1<\varepsilon$. So set
\[
\varepsilon(F_0, F_1)=\max\{\varepsilon:\: F_0, F_1\ \text{ }\ B(F_0,F_1),\: \text{}\
B(F_1,F_0)\ \text{}\ \varepsilon\}.
\]
Let $\Theta_\varepsilon(F_0)$ be the class of the continuous
distribution functions $F_1$ satisfying either $B(F_1,F_0)$ or
$B(F_0,F_1)$ with $\varepsilon(F_0, F_1)\ge\varepsilon$ and consider
another alternative hypothesis $H_{1, \varepsilon}\colon
F\in\Theta_\varepsilon(F_0)$.

Let us show, that the conditional distribution of the statistic
$R_{k,n}$ differs depending on whether $H_0$ or $H_1$ holds, that
make it possible to propose the discrimination test of these
hypotheses. The following results discusses the asymptotic behavior
of the statistic $R_{k,n}$ as $k,n\to\infty$, $k<n$, if $H_0$ holds.

\begin{theorem}\label{Th1}
If $H_0$ holds, then
\[
\sqrt{k} (R_{k,n}-1) \xrightarrow{d\,}\xi\quad \text{as}\ k, n\to\infty,
\]
\/ where $\xi$ is standard normal, $\xi\sim\operatorname{N}(0,1)$.
\end{theorem}

The above theorem allows us to propose the test distinguishing
hypotheses $H_0$ and $H_1$ for the tail of the d.f. $F$:
\begin{equation}\label{test2}
\text{if}\ R_{k,n}\notin\left(1+\frac{u_{\alpha/2}}{\sqrt{k}},
1+\frac{u_{1-\alpha/2}}{\sqrt{k}}\right),\: \text{then }\ H_0\
\text{is rejected},
\end{equation}
where $u_{\alpha/2}$ and $u_{1-\alpha/2}$ are quantiles of the
corresponding levels of the standard normal distribution. It is easy
to see, that the test is asymptotical with the significance level
$\alpha$.

Next, the following result shows the consistency of the proposed
test. Suppose that $x_{F_0}^*=x_{F_1}^*=+\infty\ \forall
F_1\in\Theta(F_0)$ (how to distinguish distributions with
$x_{F_0}^*\ne x_{F_1}^*$, see \cite{Falk,Dehaan}). Suppose in
addition, that $H_0$ does not hold and the tail of the d.f. $F$
coincides with the tail of some d.f. $F_1$, but not with the tail of
$F_0$. The consistency of the test \eqref{test2} is shown in

\begin{theorem}\label{Th2}
\begin{enumerate}[\rm(i)]
\item\vskip-5pt
If $H_1$ holds, then
\[
\sqrt{k}\,|R_{k,n}-1| \xrightarrow{d\,} +\infty
\]
as $k=k(n)\to\infty$, $k/n\to0,$ $n\to\infty$.
\item
If $H_{1,\varepsilon}$ holds, then under the same conditions
\[
\inf\limits_{F_1\in\Theta_{\varepsilon}(F_0)}\sqrt{k}\,|R_{k,n}-1| \xrightarrow{d\,}
+\infty.
\]
\end{enumerate}
\end{theorem}

Theorems~\ref{Th1}, \ref{Th2} justify the correctness of the
statistical procedure proposed in \eqref{test}. Consider two classes
$A_0$ and $A_1$ of distribution tails. Consider the null hypothesis
$\tilde{H}_0\colon F\in A_0$ and the alternative $\tilde{H}_1\colon
F\in A_1$. Suppose that there exist such d.f. $F_0$, that the
conditions $B(G, F_0)$ and $B(F_0, H)$ are satisfied with some (may
be different) $\varepsilon$ and $x_0$ for arbitrary distribution
functions $G\in A_0$ and $H\in A_1$.

\begin{corollary}\label{C1}
.\/ The procedure \eqref{test} testing the null hypothesis
$\tilde{H}_0$ \/has the asymptotical significance level $\alpha$.
Moreover, the test is consistent if $\tilde{H}_1$ holds.
\end{corollary}

Let us return to discussing the test \eqref{test2}. The proposed
test allows us to distinguish the tails of two normal distributions
with different variances, but we should weaken the B-condition
\eqref{cond} to be able to distinguish, for instance, the tails of
two normal distributions with the same variances an different means.
But weakening the B-condition implies imposing some conditions on
the sequence $k(n)$.

\begin{definition}
We say that distribution functions $F$ and $G$ satisfy the condition
$C(F,G)$, \/ if for some $\varepsilon>0$ and $x_0$
\begin{equation}\label{cond1}
\frac{(1-F(x))(-\ln(1-F(x)))^{\varepsilon}}{1-G(x)}\ \text{is
non-increasing as }\ x>x_0.
\end{equation}

\end{definition}
Let $\Theta'(F_0)$ be the class of the continuous distribution
functions $F_1$ satisfying either $C(F_1, F_0)$ or $C(F_0, F_1)$ and
the following condition: for some $\delta\in(0,1)$
\begin{equation}\label{new}
1-F_1(x)\le(1-F_0(x))^{\delta},\ \ x>x_0.
\end{equation}
It is easy to check that if $F_0$, $F_1$ satisfy either $C(F_0,F_1)$
or $C(F_1,F_0)$ with some $\varepsilon>0$, then the same holds for
all ~$\varepsilon_1$, $0<\varepsilon_1<\varepsilon$. Set
\[
\varepsilon'(F_0, F_1)=\max\{\varepsilon:\: F_0, F_1\ \text{ }\ C(F_0,F_1),\:
\text{}\ C(F_1,F_0) \text{}\ \varepsilon\}.
\]
Let $\Theta_\varepsilon'(F_0)$ be the class of the continuous
distribution functions $F_1$ satisfying \eqref{new} and either
$C(F_1,F_0)$ or $C(F_0,F_1)$ with $\varepsilon'(F_0,
F_1)\ge\varepsilon$. As before, consider the simple hypothesis
$H_0\colon F=F_0$ and two alternatives $H_1'\colon
F\in\Theta'(F_0)$, $H_{1,\varepsilon}'\colon
F\in\Theta'_\varepsilon(F_0)$, suppose in addition that $F_0$ is
continuous.

\begin{theorem}\label{Th3}
\begin{enumerate}[\rm(i)]
\item\vskip-5pt
If $H_1'$ holds, then,
\[
\sqrt{k}\,|R_{k,n}-1| \xrightarrow{d\,} +\infty
\]
as $k/n\to0$, $k^{1/2-\alpha}/\ln n\to+\infty$ \/ for some
$\alpha\in(0, 1/2),$ $n\to\infty$.
\item
If $H_{1,\varepsilon}'$ holds, then under the same conditions
\[
\inf_{F_1\in\Theta'_\varepsilon(F_0)}\sqrt{k}\,|R_{k,n}-1| \xrightarrow{d\,}
+\infty.
\]
\end{enumerate}
\end{theorem}

\section{Auxiliary results and proofs }

\subsection{Auxiliary results}

Since the statistic $R_{k,n}$ depends only on higher order
statistics of the sample, we cannot use directly the independence of
the random variables $(X_1,\ldots,X_n)$. So we consider the
conditional distribution of the statistics $R_{k,n}$ given
$X_{(n-k)}=q$ with the help of the following result.

{\renewcommand{\thelemma}{\rm\cite[~3.4.1]{Dehaan}}
\begin{lemma} Let
$X,X_1,\ldots,X_n$ be i.i.d. random variables with a common
distribution function $F$ and let $X_{(1)}\le\ldots\le X_{(n)}$ be
the $n$th order statistics. The joint distribution of the set of
statistics $\left\{\smash{X_{(i)}}\right\}_{i=n-k+1}^n$ given
$X_{(n-k)}=q$ \/ for some $k=1,\ldots,n-1$ \/  \/ agrees with the
joint distribution of the set of order statistics
$\left\{\smash{X_{(i)}^*}\right\}_{i=1}^k$ \/ of i.i.d. random
variables $\left\{X_i^*\right\}_{i=1}^k$ with the common d.f.
\[
F_q(x)=P(X\le x| X>q)=\frac{F(x)-F(q)}{1-F(q)},\quad x>q.
\]
\end{lemma}}

Let us call $F_q(x)$, $x>q$, the tail distribution function
associated with the d.f. $F$. Consider two d.f. $F$, $G$ and the
random variable $\xi_q$ with the d.f. $G_q$, $q\in\mathbb{R}$. We
set
\[
\eta_q=\ln\left(\frac{1-F(q)}{1-F(\xi_q)}\right).
\]
It is easy to check, that $\eta_q\ge0$ for all $q\in\mathbb{R}$. The
crucial point in the proof of Theorem \ref{Th2} is an investigation
of the behavior of the random variable $\eta_q$.

\begin{proposition}
Let $F_q,$ $G_q$ be the tail distribution functions associated with
$F$ and $G$ respectively. Then
\begin{enumerate}
\item
If $F_q(x)=G_q(x)$ for some $x_0$, $q>x_0$, and all $x>q$, then
$\eta_q$ is standard exponential.
\item
$G_q(x)\ge F_q(x)$ for all $x>q$ if and only if \/ $\eta_q$ is
stochastically smaller than the standard exponential random
variable.
\newline $G_q(x)\le F_q(x)$ for all $x>q$ if and only if
\/ $\eta_q$ is stochastically greater than the standard exponential
random variable.
\item
$G_q(x)\ge F_q(x)$ for all $x>q\ge x_0$ and some $x_0$ if and only
if $(1-G(x))/(1-F(x))$ is the nonincreasing function as $x>x_0$.
\end{enumerate}
\end{proposition}

\subsection{The proof of the proposition.}

(i) Suppose that $F_q(x)=G_q(x)$ for all $x>q$, then for the d.f. of
the random variable $\eta_q$ we have
\begin{equation}\label{Petaq}
\begin{aligned}[b]
P(\eta_q\le y)&=P\left(\ln\left(\frac{1-F(q)}{1-F(\xi_q)}\right)\le y\right)
=P\left(\frac{1-F(q)}{1-F(\xi_q)}\le e^y\right)\\ &=
P\left(F(\xi_q)\le1-(1-F(q))e^{-y}\right) =P\left(\xi_q\le
F^{\leftarrow}\left(1-\frac{1-F(q)}{e^y}\right)\right),
\end{aligned}
\end{equation}
with $F^{\leftarrow}(x)=\inf\{z:\: F(z)=x\}$. Besides, for the same
value of $y$
\[
P\left(\xi_q\le
F^{\leftarrow}\left(1-\frac{1-F(q)}{e^y}\right)\right)
=\frac{F\Bigl(F^{\leftarrow}\Bigl(1-\frac{1-F(q)}{e^y}\Bigr)\Bigr)
-F(q)}{1-F(q)}=1-e^{-y}.
\]

(ii) Suppose that for all $x>q$ and some $q\in\mathbb{R}$ the
relation $G_q(x)\ge F_q(x)$ holds. From \eqref{Petaq} and the
relation $1-(1-F(q))e^{-y}\ge F(q)$ with $y\ge0$ it follows that
\begin{equation}\label{F<G}
\begin{aligned}[b]
P(\eta_q\le
y)&=\frac{G\Bigl(F^{\leftarrow}\Bigl(1-\frac{1-F(q)}{e^y}\Bigr)\Bigr)
-G(q)}{1-G(q)}\\ &\ge \frac{F\Bigl(F^{\leftarrow}
\Bigl(1-\frac{1-F(q)}{e^y}\Bigr)\Bigr)-F(q)}{1-F(q)}=1-e^{-y},
\end{aligned}
\end{equation}
the result required. Let us prove the assertion in the opposite
direction. Suppose that $\eta_q$ is stochastically smaller than the
standard exponential random variable, i.e. $P(\eta_q\le
y)\ge1-e^{-y}$ with $y\ge0$. Using \eqref{Petaq}, we have

\begin{align*}
\frac{G\Bigl(F^{\leftarrow}\Bigl(1-\frac{1-F(q)}{e^y}\Bigr)\Bigr)
-G(q)}{1-G(q)}\ge1-e^{-y} &\Longleftrightarrow
\frac{1-G\Bigl(F^{\leftarrow}\Bigl(1-\frac{1-F(q)}{e^y}\Bigr)\Bigr)}{1-G(q)}\le
e^{-y} \\ &\Longleftrightarrow
G\left(F^{\leftarrow}\left(1-\frac{1-F(q)}{e^y}\right)
\right)\le1-\frac{1-G(q)}{e^y}\\ &\Longleftrightarrow
F^{\leftarrow}\left(1-\frac{1-F(q)}{e^y}\right)\le
G^{\leftarrow}\left(1-\frac {1-G(q)}{e^y}\right).
\end{align*}
Denote $z_{F}=F^{\leftarrow}\left(1-e^{-y}(1-F(q))\right)$ and
$z_{G}=G^{\leftarrow}\left(1-e^{-y}(1-G(q))\right)$. Since
$F(z_{F})=1-e^{-y}(1-F(q))$ and $G(z_{G})=1-e^{-y}(1-G(q))$, we have
\[
e^{-y}=\frac{1-G(z_{G})}{1-G(q)}=\frac{1-F(z_{F})}{1-F(q)}.
\]
Next, since $z_{F}\le z_{G}$,
\[
\frac{1-F(z_{F})}{1-F(q)}=\frac{1-G(z_{G})}{1-G(q)}\le\frac{1-G(z_{F})}{1-G(q)}.
\]
This observation ends the proof, since $z_{F}\in\lbrack q,\infty)$.
The proof of the second assertion of the item (ii) is similar.

(iii) We have
\begin{align*}
\frac{G(x)-G(q)}{1-G(q)}\ge\frac{F(x)-F(q)}{1-F(q)}\quad \forall
x>q\ge x_0
&\Longleftrightarrow\frac{1-G(x)}{1-G(q)}\le\frac{1-F(x)}{1-F(q)}\quad
\forall x>q\ge x_0\\ &\Longleftrightarrow
\frac{1-G(x)}{1-F(x)}\le\frac{1-G(q)}{1-F(q)}\quad \forall x>q\ge
x_0\\ &\Longleftrightarrow\frac{1-G(x)}{1-F(x)}\ \text{is
nonincreasing for all } x>x_0.
\end{align*}

\subsection{The proof of Theorem \ref{Th1}}

The proof of Theorem \ref{Th1} is evident, the same steps are used
in extreme value theory to prove the consistency of the Hill's
estimator of extreme value index (see \cite[~lemma 3.2.3]{Dehaan}
and its proof). Though, let us give the proof of Theorem \ref{Th1}.
In assumptions of Theorem \ref{Th1} $F_0(X_1)$ has a uniform
distribution on the interval $[0,1]$, $F_0(X_1)\sim U[0,1]$, hence
$-\ln(1-F_0(X))$ is standard exponential. From R\'enyi's
representation (see \cite{Dehaan}) it follows that
\[
\bigl\{-\ln(1-F_0(X_{(n-i)}))+\ln(1-F_0(X_{(n-k)}))\bigr\}_{i=0}^{k-1}
\overset{d}{=}\Biggl\{\sum\limits_{j=i+1}^k \frac{E_{n-j+1}}{j}\Biggr\}_{i=0}^{k-1},
\]
where $E_1,E_{2},\ldots\strut$ are independent standard exponential
random variables. Therefore, the distribution of the left part does
not depend on $n$ and
\[
\bigl\{-\ln(1-F_0(X_{(n-i)}))+\ln(1-F_0(X_{(n-k)}))\bigr\}_{i=0}^{k-1}
\overset{d}{=}\bigl\{E_{(k-i)}\bigr\}_{i=0}^{k-1},
\]
where $E_{(1)}\le\ldots\le E_{(k)}$ are the $k$th order statistics
of $\left\{E_i\right\}_{i=1}^k$.
\[
\sqrt{k}(R_{k,n}-1)\overset{d}{=}\sqrt{k}\left(\frac{1}{k}\sum
\limits_{i=0}^{k-1}E_{(k-i)}-1\right)
=\sqrt{k}\Biggl(\frac{1}{k}\sum_{j=1}^kE_j-1\Biggr),
\]
and the result follows from the fact that the characteristic function of the right-hand side tends to the characteristic function of the standard normal distribution as $k\to\infty.$

\subsection{The proof of Theorem \ref{Th2}}

Let us prove the item (i) at first. The scheme of the proof includes
the steps that are similar to ones using in works
\cite{Rodionov1,Rodionov2}. Consider the asymptotic behavour of the
statistic $R_{k,n}$ as $n\to\infty$. Denote
\[
Y_i=\ln(1-F_0(q))-\ln(1-F_0(X_i^*)),
\] where
$\left\{X_i^*\right\}_{i=1}^k$ are i.i.d. random variables
introduced in Lemma with the d.f.
\[
F_q(x)=\frac{F_1(x)-F_1(q)}{1-F_1(q)},\quad q<x.
\]
Assuming $F=F_0$ and $G=F_1$, we have $Y_i\overset{d}{=}\eta_q$,
$i\in\{1,\ldots,k\}$. It follows from Lemma, that the joint
conditional distribution of the $k$th order statistics
$\left\{\smash{Y_{(i)}}\right\}_{i=1}^k$ of the sample
$\left\{Y_j\right\}_{i=1}^k$ given $X_{(n-k)}=q$ equals the joint
distribution of the $k$th order statistics
$\left\{\smash{Z_{(j)}}\right\}_{j=1}^k$ of the sample
$\left\{\smash{Z_j}\right\}_{j=1}^k$, where
\[
Z_j=\ln(1-F_0(X_{(n-k)}))-\ln(1-F_0(X_{(n-j+1)})),\quad j=1,\ldots,k.
\]
Next, we see that
\[
R_{k,n}=\frac{1}{k}\sum_{i=1}^kZ_i.
\]
Therefore, the conditional distribution of $R_{k,n}$ given
$X_{(n-k)}=q$ agrees with the distribution of the statistic
$\frac{1}{k}\sum\limits_{i=1}^kY_i$. Next, the distribution
functions $F_1$ and $F_0$ satisfy either the condition $B(F_0,F_1)$
or the condition $B(F_1,F_0)$ under the assumptions of Theorem 2.
Firstly suppose that $B(F_0,F_1)$ is satisfied for some
$\varepsilon>0$ and $x_0$. Since $x^*=+\infty$, then
$X_{(n-k_n)}\to+\infty$ a.s., so we can only consider the case
$q>x_0$. The item (iii) of Proposition implies, that
\[
\frac{1-F_1(x)}{1-F_1(x_0)}\ge
\frac{(1-F_0(x))^{1-\varepsilon}}{(1-F_0(x_0))^{1-\varepsilon}},\quad x>x_0.
\]
Using \eqref{F<G}, we have
\begin{align*}
P(Y_1\le x)&=1-\frac{1-F_1\Bigl(F_0^{\leftarrow}\Bigl(1-\frac
{1-F_0(q)}{e^x}\Bigr)\Bigr)}{1-F_1(q)}\\ &\le
1-\frac{\Bigl(1-F_0\Bigl(F_0^{\leftarrow}\Bigl(1-\frac{1-F_0(q)}{e^x}\Bigr)\Bigr)
\Bigr)^{1-\varepsilon}}{(1-F_0(q))^{1-\varepsilon}}=1-e^{-(1-\varepsilon)x},
\end{align*}
therefore $Y_1$ is stochastically greater than a random variable
$E\sim \Exp(1-\varepsilon)$, written $Y_1\gg E$. Next, let
$E_1,\ldots,E_{k_n}$ be i.i.d. random variables with the common d.f.
$H(x)=1-e^{-(1-\varepsilon)x}$, then
\begin{equation}\label{ll}
\sqrt{k}\left(\frac{1}{k}\sum_{i=1}^kY_i-1\right) \gg
\sqrt{k}\left(\frac{1}{k}\sum_{i=1}^kE_i-1\right).
\end{equation}
Since \eqref{ll} holds for all $q>x_0$ and $X_{(n-k)}\to +\infty$
a.s. as $n\to\infty$, then under the conditions of Theorem \ref{Th2}
we have
\begin{equation}\label{LL}
\sqrt{k}(R_{k,n}-1)\gg\sqrt{k}\left(\frac{1}{k}\sum_{i=1}^kE_i-1\right).
\end{equation}
From the Lindeberg--Feller central limit theorem,
\[
(1-\varepsilon)\sqrt{k}\left(\frac{1}{k}\sum_{i=1}^k
E_i-\frac{1}{1-\varepsilon}\right) \xrightarrow{d}\xi\sim\operatorname{N}(0,1),\quad
n\to\infty,
\] hence

\begin{equation}\label{ref2}
\sqrt{k}\left(\frac{1}{k}\sum_{i=1}^kE_i-1\right) \xrightarrow{P}+\infty,\quad
n\to\infty.
\end{equation}
Finally, using \eqref{LL}, we have
\begin{equation}\label{C2}
\sqrt{k}(R_{k,n}-1)\xrightarrow{P}+\infty,\quad n\to\infty.
\end{equation}
If the condition $B(F_1,F_0)$ holds, then using the similar steps,
we have
\begin{equation}\label{C3}
\sqrt{k}(R_{k,n}-1)\xrightarrow{P}-\infty,\quad n\to\infty.
\end{equation}
The second assertion of the theorem follows easily from \eqref{LL}
and \eqref{ref2}.

\subsection{The proof of the corollary}

As it follows from \eqref{C2} and \eqref{C3}, if $\tilde{H}_0$
holds, then $\sqrt{k}(R_{k,n}-1)\xrightarrow{P}+\infty$ as
$n\to\infty$. On the other hand, if $\tilde{H}_1$ holds, then
$\sqrt{k}(R_{k,n}-1){\xrightarrow{P}}-\infty$ as $n\to\infty$, and
it proves that the test \eqref{test} is consistent.

\subsection{The proof of Theorem \ref{Th3}}

Let us prove the item (i) at first. Denote $\overline{F}(x)=1-F(x)$.
Using the notation of Theorem \ref{Th2}, find the distribution of
the random variable $Y_1$. Suppose, that the condition $C(F_0,F_1)$
is satisfied. Using \eqref{F<G} and the item (iii) of proposition,
\begin{align*}
P(Y_1\le x)&=1-\frac{\Fbarl\left({\FbarO}^{\leftarrow}
(\FbarO(q){e^{-x}})\right)}{\Fbarl(q)}\\ &\le
1-\frac{\FbarO(q)e^{-x}\left(-\ln(\FbarO(q)e^{-x})\right)
^{\varepsilon}}{\overline{F}(q)(-\ln\FbarO(q))^{\varepsilon}}
=1-e^{-x}\left(1+\frac{x}{-\ln\FbarO(q)}\right)^{\varepsilon}.
\end{align*}
For $\varepsilon,c\in(0,1)$ we have
\[
(1+cx)^{\varepsilon}\ge1+c\varepsilon-c\varepsilon e^{-x},\quad x\ge0,
\]
and $G(x)=1-e^{-x}(1+c\varepsilon-c\varepsilon e^{-x})$ is the
distribution function. Therefore,
\[
P(Y_1\le x)\ge1-e^{-x}-\frac{\varepsilon}{-\ln\FbarO(q)}(1-e^{-x}).
\]
Further, let $\zeta,\zeta_1,\ldots,\zeta_k$ be i.i.d. random
variables with such d.f. Similarly to the proof of Theorem
\ref{Th2}, we have
\begin{equation}\label{ll1}
\sqrt{k}\left(\frac{1}{k}\sum_{i=1}^kY_i-1\right) \gg
\sqrt{k}\left(\frac{1}{k}\sum_{i=1}^k\zeta_i-1\right).
\end{equation}
It is easy to see, that
\[
E\zeta=1+\frac{\varepsilon}{-2\ln \overline{F}_0(q)},
\operatorname{Var}\zeta=1-\left(\frac{\varepsilon}{-2\ln\overline{F}_0(q)}\right)^2,
\]
so we obtain
\begin{equation}\label{ref3}
\sqrt{k}\left(\frac{1}{k}\sum_{i=1}^k\zeta_i-1\right)
=\sqrt{k}\left(\frac{1}{k}\sum_{i=1}^k\zeta_i-E\zeta\right)
+\sqrt{k}\frac{\varepsilon}{-2\ln\overline{F}_0(q)}.
\end{equation}

Now consider the statistic $\sqrt{k}/\ln\overline{F}_0(X_{(n-k)})$.
Introduce $R_i=\overline{F}_1(X_i)$, $i=1,\ldots,n$. Since $F_1$ is
continuous, then $R_1,\ldots,R_n$ are independent standard uniform
random variables and $R_{(k)}=\overline{F}_1(X_{(n-k)})$. An appeal
to \cite[~theorem 2.2.1]{Dehaan} shows that
\begin{equation}\label{last}
\frac{n}{\sqrt{k}}\left(R_{(k)}-\frac{k}{n}\right)
\xrightarrow{d\,}\operatorname{N}(0,1),\quad n\to\infty.
\end{equation}
Applying the Delta Method (see, for instance, \cite{Resnick}) with
the function $f(x)=-x/\ln x$, we have
\[
\frac{n}{\sqrt{k}}\left(\frac{R_{(k)}}{-\ln R_{(k)}}-\frac {k/n}{-\ln(k/n)}\right)
\xrightarrow{P}0,\ n\to\infty,
\]
therefore under the conditions of Theorem 3,
\[
f'\biggl(\frac{k}{n}\biggr) =-\frac{1}{\ln(n/k)}+\frac {1}{(\ln(n/k))^2}\to0,\quad
n\to\infty.
\]
Further,
\[
\frac{n}{\sqrt{k}}\left(\frac{R_{(k)}}{\ln R_{(k)}}-\frac
{k/n}{\ln(k/n)}\right)=\frac{n}{\sqrt{k}}\left(\frac {R_{(k)}}{\ln
R_{(k)}}-\frac{k/n}{\ln(R_{(k)})}\right) +\sqrt{k}\left(\frac{1}{\ln
R_{(k)}}-\frac{1}{\ln(k/n)}\right),
\]
and from \eqref{last} it follows that the first summand in the right
hand converges to $0$ in probability as $n\to\infty$. Thus,
\[
\sqrt{k}\left(\frac{1}{\ln R_{(k)}}-\frac{1}{\ln(k/n)}\right) \xrightarrow{P}0,\quad
n\to\infty,
\]
and under conditions of Theorem \ref{Th3}
\[
\frac{\sqrt{k}}{-\ln R_{(k)}}=\sqrt{k}\left(\frac{1}{-\ln
R_{(k)}}-\frac{1}{-\ln(k/n)}\right) +\frac{\sqrt{k}}{-\ln
(k/n)}\xrightarrow{P}+\infty,\quad n\to\infty.
\]
On the other hand, using \eqref{new}, we have
\begin{equation}\label{ref4}
\frac{\sqrt{k}}{-\ln\overline{F}_0(X_{(n-k)})}
=\frac{\sqrt{k}}{-\ln\overline{F}_0\left({\overline{F}_1}^{\leftarrow}(R_{(k)})\right)}
\ge\frac{\sqrt{k}}{-\delta^{-1}\ln R_{(k)}}\xrightarrow{P}+\infty
\end{equation}
as $n\to\infty$. Next, from the Law of Large Numbers for triangular
arrays (see \cite{Mik}) it follows that for all $\epsilon>0$
\[
\sqrt{k}\left(\frac{1}{k}\sum_{i=1}^k\zeta_i-1\right) =o_P^{}(k^{\epsilon}),\quad
n\to\infty,
\]
i.e. the term in the left hand is asymptotically smaller in
probability than $k^{\epsilon}$. Thus, for all $q$ given
$X_{(n-k)}=q$
\[
\sqrt{k}\left(\frac{1}{k}\sum_{i=1}^kY_i-1\right) \xrightarrow{P}+\infty,\quad
n\to\infty,
\]
and finally
\[
\sqrt{k}\left(R_{k,n}-1\right) \xrightarrow{P}+\infty,\quad n\to\infty.
\]
If the condition $C(F_1,F_0)$ is satisfied, then
\[
\sqrt{k_n}(R_{k,n}-1)\xrightarrow{P}-\infty,\quad n\to\infty,
\]
and the proof is similar. The second assertion of Theorem follows
easily from \eqref{ll1}, \eqref{ref3} and \eqref{ref4}.

\smallskip
The author is grateful to Prof. V.~I. Piterbarg for his invaluable
help in the preparation of the manuscript.

\end{document}